\documentclass[11pt, reqno]{amsart}
\usepackage{amsmath}
\usepackage{amssymb}
\input amssym.def
\input amssym.tex

\font\tengothic=eufm10
\font\sevengothic=eufm7
\newfam\gothicfam
      \textfont\gothicfam=\tengothic
      \scriptfont\gothicfam=\sevengothic

\numberwithin{equation}{section}

\begin{document}
\setlength{\baselineskip}{1.3em}

{\theoremstyle{plain}
    \newtheorem{thm}{\bf Theorem}[section]
    \newtheorem{pro}[thm]{\bf Proposition}
    \newtheorem{lemma}[thm]{\bf Lemma}
    \newtheorem{cor}[thm]{\bf Corollary}
    \newtheorem{claim}[thm]{\bf Claim}
}
{\theoremstyle{remark}
    \newtheorem{conj}[thm]{\bf Conjecture}
    \newtheorem{remark}[thm]{\bf Remark}
    \newtheorem{problem}[thm]{\bf Problem}
}
{\theoremstyle{definition}
    \newtheorem{defn}[thm]{\bf Definition}
}

\newcommand{\cc}{{\mathbb C}}
\newcommand{\pp}{{\mathbb P}}
\newcommand{\zz}{{\mathbb Z}}
\newcommand{\nn}{{\mathbb N}}
\newcommand{\R}{{\mathcal R}}
\newcommand{\A}{{\mathcal A}}
\newcommand{\x}{{\mathbb X}}
\newcommand{\y}{{\mathbb Y}}
\newcommand{\ix}{I_{\mathbb X}}
\newcommand{\I}{{\mathcal I}}
\newcommand{\B}{{\bf B}}
\newcommand{\bi}{{\bf I}}
\newcommand{\V}{{\bf V}}
\newcommand{\cv}{{\mathcal V}}
\def\L{\mathcal L}
\def\O{\mathcal O}
\def\F{\mathcal F}
\def\K{\mathcal K}
\def\M{\mathcal M}
\def\N{\mathcal N}
\def\J{\mathcal I}
\def\H{{\bf H}}
\def\m{\frak m}
\def\To{\longrightarrow}
\def\reg{\operatorname{reg}}
\newcommand{\spec}{\mbox{\rm Spec }}
\newcommand{\proj}{\mbox{\rm Proj }}
\newcommand{\Hom}{\mbox{\rm Hom}}
\def\tor{\operatorname{Tor}}
\newcommand{\ext}{\mbox{\rm Ext}}
\newcommand{\smap}{\longrightarrow\!\!\!\!\!\!\!\!\!\longrightarrow}
\newcommand{\sfrac}[2]{\frac{\displaystyle #1}{\displaystyle #2}}
\newcommand{\under}[1]{\underline{#1}}
\newcommand{\ov}[1]{\overline{#1}}
\newcommand{\sheaf}[1]{{\mathcal #1}}

\title{Adjoint line bundles and syzygies of projective varieties}
\author{Huy T\`ai H\`a}
\address{Tulane University, Department of Mathematics, 6823 St. Charles Ave., New Orleans, LA 70118}
\email{tai@math.tulane.edu}

\begin{abstract} Let $X$ be a smooth projective variety and let $K$ be the canonical divisor of $X$. In this paper, we study embeddings of $X$ given by adjoint line bundles $K \otimes L$, where $L$ is an ample invertible sheaf. When $X$ is a regular surface and $K \otimes L$ is normally generated (i.e. satisfies property $N_0$), we obtain a numerical criterion for $K \otimes L$ to have property $N_p$. When $X$ is a regular variety of arbitrary dimension, under a mild condition, we give an explicit calculation for the regularity of ideal sheaves of such embeddings.
\end{abstract}
\maketitle


\begin{center}
{\it Dedicated to Prof. H\`a Huy Kho\'ai on the occasion of his sixtieth birthday}
\end{center}

\section{Introduction}

\noindent A classical result of Castelnuovo showed that a curve of degree greater than $2g$ has a normal homogeneous coordinate ring ($g$ is the genus of the curve). This result was rediscovered and further extended by many authors including Mumford \cite{mum}, St. Donat \cite{std} and Green \cite{gr}. Making use of techniques of Koszul cohomology, Green \cite[Theorem 4.a.1]{gr} proved that a curve of degree greater than $(2g+p)$ has a normal homogeneous coordinate ring, and its defining ideal is generated by quadrics and has linear syzygies up to the $p$-th step. This property was later developed and called {\it property $N_p$} \cite{gl,el}.

Mukai observed that the studies on curves can be viewed as dealing with {\it adjoint bundles} of the type $K \otimes D$, where $K$ is the canonical divisor and $D$ is an ample line bundle. In particular, he conjectured that if $X$ is a smooth projective surface with canonical divisor $K_X$, and $L$ is an ample line bundle on $X$, then for any integer $p \ge 0$, the adjoint bundle $K_X \otimes L^{\otimes (4+p)}$ has property $N_p$. In higher dimension, Ein and Lazarsfeld \cite{el} extends Mukai's conjecture as follows: if $X$ is a smooth complex projective variety of dimension $n$ with canonical divisor $K_X$, and $L$ is an ample line bundle on $X$, then for any integer $p \ge 0$, the adjoint bundle $K_X \otimes L^{\otimes (n+2+p)}$ has property $N_p$. For abelian varieties, Lazarsfeld \cite{la1} conjectured that if $X$ is an abelian variety and $L$ is an ample line bundle on $X$, then for any integer $p \ge 0$, the adjoint bundle $K_X \otimes L^{\otimes (p+3)}$ has property $N_p$.  

A considerable amount of efforts has since then been put forward to investigate the minimal free resolution of embeddings of projective varieties given by adjoint line bundles. Particularly, among other things, Gallego and Purnaprajna in a series of papers \cite{gp1, gp2, gp3, gp4, gp5, gp6} give many special situations in which Mukai's conjecture holds true; Park \cite{p1, p2} extended many of Gallego and Purnaprajna's results for ruled surfaces; Ein and Lazarsfeld \cite{el} prove their conjecture when the ample line bundle $L$ is very ample; Hering, Schenck and Smith \cite{ss} establish Ein and Lazarsfeld's conjecture for Gorenstein projective toric varieties; Rubei \cite{rubei}, Pareschi \cite{pa}, and Pareschi and Popa \cite{pp1, pp2} study syzygies of adjoint line bundles on abelian varieties in connection to Lazarsfeld's conjecture. 

In this paper, we study adjoint line bundles on {\it regular} projective varieties. We refer the reader to Definition \ref{aCM-def} for the notion of regular varieties. Let $X$ be a smooth regular projective variety, and let $K$ and $L$ be the canonical divisor and an ample line bundle on $X$. Let $X_\L$ be the projective embedding of $X$ given by the adjoint line bundle $\L = K \otimes L$ (if $\L$ is very ample), and let $\J_\L$ be the ideal sheaf of $X_\L$. Similar to previous work, we investigate the minimal free resolution of $X_\L$; and thus, in particular, we examine property $N_p$ of $\L$. We are also interested in the Castelnuovo-Mumford regularity of $\J_\L$.

Section \ref{section-np} is devoted to the case when $X$ is a regular surface. Our work in this section is inspired by a result of Gallego and Purnaprajna. \cite[Theorem 1.3]{gp6} states that if $X$ is a rational surface, $L$ is a globally generated ample line bundle on $X$, then $L$ has property $N_p$ provided $-K.L \ge p+3$. For adjoint line bundles, this shows that if $L$ is an ample line bundle on a rational surface $X$ such that $K \otimes L$ is globally generated, then $K \otimes L$ has property $N_p$ if $-K.L \ge p+3+K^2$. We extend this result to regular surfaces, and seek for an improvement on the bound $p+3+K^2$. Let $X$ be a regular surface, and let $p_a$ be the arithmetic genus of $X$. Then, $p_a = \chi(\O_X) - 1 = h^2(X, \O_X) \ge 0$. Our main result, Theorem \ref{surface}, shows that the bound $-K.L \ge p+3+K^2$ can be strengthened to be $-K.L \ge p+3+K^2 - 2p_a$. It was pointed out to the author by Purnaprajna that, under our condition on $-K.L$, regular anticanonical surfaces are rational. As a consequence of Theorem \ref{surface}, we prove Mukai's conjecture for rational anticanonical surfaces. This is done in Theorem \ref{surface-np} and Corollary \ref{surface-cor}. Theorem \ref{surface-np} provides a different proof for the first part of \cite[Theorem 1.23]{gp6}. Purnaprajna's observation also simplifies and fixes a mistake in our original proof of Lemma \ref{normal}. 

In Section \ref{section-reg}, we study the Castelnuovo-Mumford regularity of $\J_\L$. Let $n = \dim X$, we show that when the adjoint line bundle $\L$ is normally generated, under a mild condition, $\reg \J_\L$ can be explicitly calculated, namely $\reg \J_\L = n+2$ if $H^0(X, K) \not= 0$ and $\reg \J_\L = n+1$ if $H^0(X, K) = 0$. This is done in Theorem \ref{aCM}. Our theorem can be interpreted as a corresponding version for adjoint line bundles of Green's Top Row theorem \cite[Theorem 4.a.4]{gr}. This result also exhibits an interesting asymptotic behaviour of adjoint line bundles of the form $K \otimes L^{\otimes m}$. As Remark \ref{behaviour} indicates, when $m$ gets large, while the beginning part of the minimal free resolution of $\J_\L$ becomes more linear due to Ein and Lazarsfeld's theorem, its width is always fixed; in other words, even though the beginning part of the minimal free resolution of $\J_\L$ becomes nicer, the level of complexity of $\J_\L$ remains the same.

Our approach in this paper is based on the classical method of investigating general hyperplane sections. This idea grows out of our previous work \cite{ha1}. Theorem \ref{surface} is proved by closely examining the Hilbert function and Betti numbers of general hyperplane sections of curves and surfaces. Theorem \ref{surface-np} is attained by combining Theorem \ref{surface} (and in particular, Corollary \ref{surface-main}) and the argument of Gallego and Purnaprajna \cite[Proposition 1.10]{gp6}. Theorem \ref{aCM} is established by combining vanishing theorems with the study of Koszul cohomology groups, making use of Green's work in \cite{gr} to derive the non-vanishing of certain Koszul cohomology groups.

\noindent{\it Acknowledgement.} The author would like to thank B.P. Purnaprajna for the observation about regular anticanonical and rational anticanonical surfaces, which in particular fixes a mistake in our proof of Lemma \ref{normal} (which appeared in the first draft of the paper).


\section{Notations and terminology} \label{section-intro} 

\noindent Our notations and terminology mostly follow from \cite{hart, l}. Throughout the paper, we will work over an algebraically closed field $k$ of characteristic 0. Although most of our results hold for unmixed smooth projective schemes, we shall only discuss varieties. By a {\it projective variety} we shall mean a reduced irreducible smooth projective scheme over $k$. A {\it surface} is a projective variety of dimension 2. To avoid the tension between the {\it additive language} of divisors and the {\it multiplicative formalism} of line bundles, our convention is always to work multiplicatively. And so, in many occasions, we identify divisors with corresponding line bundles. For example, if $X$ is a smooth projective variety with structure sheaf $\O_X$ and canonical divisor $K$, and suppose $L$ is an ample line bundle on $X$, then by writing $K \otimes L$ we shall mean $\omega \otimes L$ where $\omega = \O_X(K)$ is the dualizing sheaf of $X$. On the other hand, we sometimes also use $K + L$ to refer to $K + D$, where $D$ is a divisor such that $L = \O_X(D)$. Note that throughout the paper, the additive notation (such as $K+L$) is used only to calculate intersection numbers which are invariant under numerical equivalence; and thus, no confusion will be caused.

We shall define regular projective varieties. This class of varieties includes rational projective varieties, and more generally, projective varieties having arithmetically Cohen-Macaulay embeddings.

\begin{defn} \label{aCM-def}
A projective variety $X$ of dimension $n$ is said to be {\it regular} if $$H^{n-1}(X, \O_X) = 0.$$ 
\end{defn}

By Serre's duality, a projective variety $X$ is regular if $H^1(X, K) = 0$, where $K$ is the canonical divisor of $X$.

\begin{remark} \label{rational}
A smooth projective variety $X$ of dimension $n$ is {\it rational} if it is birationally equivalent to $\pp^n$, i.e. there exists a birational morphism $\pi: X \To \pp^n$. Since $X$ is smooth, it follows from \cite[Theorem 1.2]{cha} that 
$$H^i(X, \O_X) \simeq H^i(\pp^n, \O_{\pp^n}) = 0 \ \forall \ i \ge 1.$$
Thus, a rational projective variety is regular.
\end{remark}

\begin{remark} A smooth projective variety $X$ is said to have an {\it arithmetically Cohen-Macaulay embedding} if $X$ has a projective embedding with Cohen-Macaulay coordinate ring, i.e. there exists a Cohen-Macaulay standard graded $k$-algebra $A$ such that $X \simeq \proj A$. The class of varieties having arithmetically Cohen-Macaulay embeddings has recently attracted considerably attention (cf. \cite{ggp, chtv, cha, che, ha1, ht}). Since $X$ is smooth, $X$ is a Cohen-Macaulay scheme. Thus, it follows from \cite[Lemma 1.1]{cha} that $X$ has an arithmetically Cohen-Macaulay embedding if and only if $H^i(X, \O_X) = 0$ for $i = 1, \dots, \dim X - 1$. Hence, if $X$ has an arithmetically Cohen-Macaulay embedding then $X$ is regular. 
\end{remark}

We now direct those who are familiar with the concepts of regularity, property $N_p$ and Koszul cohomology to the next section. For others, we recall these notions following \cite{mum, gr, gl, el}.

\begin{defn} Let $X$ be a projective variety, and let $\F$ be a coherent sheaf on $X$. The {\it regularity} of $\F$, denoted by $\reg \F$, is defined to be the least integer $r$ such that $H^i(X, \F(r-i)) = 0$ for all $i > 0$.
\end{defn}

The regularity can also be interpreted in terms of minimal free resolution as follows. Let $X \hookrightarrow \pp^n$ be an embedding of $X$ into a projective space, and let
$$0 \To \bigoplus_{j} \O_{\pp^n}(-b_{sj}) \To \dots \To \bigoplus_{j} \O_{\pp^n}(-b_{0j}) \To \F \To 0$$
be the minimal free resolution of $\F$ as an $\O_{\pp^n}$-module. Then 
$$\reg \F = \max \{ b_{ij} - i ~|~ b_{ij} \not= 0 \}.$$

\begin{defn} \label{np-property}
Let $X$ be a smooth projective variety, and let $\L$ be a very ample line bundle. The sections of $\L$ give an embedding $X \hookrightarrow \pp(H^0(X, \L)^*)$. Let $R = \mbox{Sym}^* H^0(X, \L)$ be the coordinate ring of $\pp(H^0(X,\L)^*)$, and let $S$ be the coordinate ring of $X$ in $\pp(H^0(X, \L)^*)$. Let 
\[ 0 \to F_s \to F_{s-1} \to \cdots \to F_0 \to S \to 0 \]
be the minimal free resolution of $S$ over $R$. 
\begin{enumerate}
\item $\L$ is said to have {\it property $N_0$} if $\L$ is normally generated, i.e. $S \simeq \oplus_{i \ge 0} H^0(X, \L^{\otimes i})$.
\item For an integer $p \ge 1$, $\L$ is said to have {\it property $N_p$} if $\L$ has property $N_0$, $F_0 = R$ and $F_i = R(-i-1)^{\alpha_i}$ for $1 \le i \le p$.
\end{enumerate}
\end{defn}

Roughly speaking, property $N_0$ says that $\L$ embeds $X$ as a projectively normal variety, property $N_1$ requires that $\L$ embeds $X$ as a projectively normal variety whose defining ideal is generated by quadrics, and more generally, for an integer $p \ge 2$, property $N_p$ means that $\L$ embeds $X$ as a projectively normal variety whose defining ideal is generated by quadrics and has linear syzygies up to the $p$-th step.

\begin{defn} \label{koszul}
Let $X$ be a projective variety. Let $\L$ be an ample line bundle and $\F$ a coherent sheaf on $X$. Let $W = H^0(X, \L)$ and $S = \mbox{Sym}^{*} W$. Then,  $S$ is the homogeneous coordinate ring of $\pp(W)$. Let $B = B(\F, \L) = \oplus_{q \in \zz} H^0(X, \F \otimes \L^{\otimes q}) = \oplus_{q \in \zz} B_q$ a $S$-graded module.
The {\it Koszul complex} of $B$ is defined to be
\[ \dots \To \bigwedge^{p+1}W \otimes B_{q-1}
 \stackrel{d_{p+1,q-1}}{\To} \bigwedge^pW \otimes B_q
\stackrel{d_{p,q}}{\To} \bigwedge^{p-1}W \otimes  B_{q+1} \To
\dots \] 
and the {\it Koszul cohomology groups} of $B$ are defined to be
\[ \K_{p,q}(\F, \L) = \sfrac{\mbox{ker } d_{p,q}}{\mbox{im } d_{p+1,q-1}},
 \mbox{ for } p,q \in \zz. \]
\end{defn}

When $\L = \O_X(D)$ is the invertible sheaf corresponding to a divisor $D$ and $\F = \widetilde{M}$ is the sheaf associated to a module $M$, we write
$\K_{p,q}(M, D)$  for $\K_{p,q}(\F, \L)$. When $\F$ is the structure sheaf $\O_X$, we also write $\K_{p,q}(\L)$ for $\K_{p,q}(\F, \L)$.


\section{Adjoint line bundles on regular surfaces} \label{section-np}

\noindent In this section, we shall derive a numerical criterion for adjoint line bundles on a regular surface to have property $N_p$, and thus, provide a new proof of Mukai's conjecture for rational anticanonical surfaces. 

The next lemma is well known. We will include the proof for self-containment.

\begin{lemma} \label{can-mod}
Let $R$ be a standard graded $k$-algebra of dimension $d$, and let $I \subset R$ be a homogeneous ideal. Suppose $A = R/I$ is of dimension $n$. Then, $H^n_\m(A)^*$ is a torsion free $A$-module, where $\m$ is the maximal homogeneous ideal of $R$.
\end{lemma} 

\begin{proof} Let $x \in A$ be a homogeneous element, which is not a zero divisor. Let $\delta$ be the degree of $x$. Consider the exact sequence of $R$-modules
$$0 \To A(-\delta) \stackrel{\times x}{\To} A \To A/xA \To 0.$$
This sequence gives rise to the following exact sequence
$$\ext^{d-n}_R(A/xA,R) \To \ext^{d-n}_R(A,R) \stackrel{\times x}{\To} \ext^{d-n}_R(A,R)(-\delta).$$
Since $\dim A/xA = n-1$, by local duality, we have
$\ext^{d-n}_R(A/xA,R) = H^n_\m(A/xA) = 0.$
Thus, we have an injection
$\ext^{d-n}_R(A,R) \hookrightarrow \ext^{d-n}_R(A,R)(-\delta).$
This implies that $\ext^{d-n}_R(A,R)$ is torsion free as $A$-module. Moreover, by local duality again, $H^n_\m(A)^* = \ext^{d-n}_R(A,R)(-d)$. The conclusion follows.
\end{proof} 

Recall from the definition that a smooth surface $X$ is regular if $H^1(X, \O_X) = 0$. Our first main result is stated as follows. Our proof is inspired by that of \cite[Theorem 2.5]{ha1}.
 
\begin{thm} \label{surface}
Let $X$ be a surface with $H^1(X, \O_X) = 0$. Suppose $K$ is the canonical divisor and $L$ is an ample divisor on $X$ such that $\L = K \otimes L$ is normally generated. Then, for any integer $p \ge 0$, the adjoint line bundle $\L$ has property $N_p$ if $-K.L \ge p+3+K^2-2p_a$, where $p_a$ is the arithmetic genus of $X$.
\end{thm}

\begin{proof} Since $\L$ is normally generated, $\L$ is necessarily very ample (cf. \cite{mum}). Thus, the theorem is proved for $p = 0$. Suppose now that $p \ge 1$. Let $N = \dim_k H^0(X, \L) - 1$, and let $X_\L$ be the embedding of $X$ in $\pp^N$ given by $\L$. Let $S, I$ and $\J$ be the coordinate ring, the defining ideal and the ideal sheaf of $X_\L$ in $\pp^N$. 

For $h > 0$, since $\L^{\otimes h} = K \otimes \L^{\otimes (h-1)} \otimes L$, by Kodaira's vanishing theorem, we have
$H^i(X, \L^{\otimes h}) = 0 \ \text{for} \ i \ge 1 \ \text{and} \ h > 0.$
This, together with Kodaira's vanishing theorem (for $h < 0$ since $\L$ is very ample), and the hypotheses (for $h = 0$), gives us
$$H^1(X, \L^{\otimes h}) = 0 \ \text{for} \ h \in \zz.$$
Thus, $\L$ embeds $X$ as an arithmetically Cohen-Macaulay variety. That is, $S/I$ is a Cohen-Macaulay ring. Let 
$$ 0 \To F_{N-2} \To \dots \To F_1 \To I \To 0$$
be the minimal free resolution of $I$, where $$F_i = \bigoplus_{j \ge i+1} S(-j)^{\beta_{i,j}} \ \text{for} \ i = 1, \dots, N-2.$$
The adjoint line bundle $\L$ has property $N_p$ if and only if $\beta_{i, j} = 0$ for all $1 \le i \le p$ and $j \ge i+2$.

Let $C$ be a general hyperplane section of $X_\L \subset \pp^N$, then $C$ is an irreducible arithmetically Cohen-Macaulay curve in $\pp^{N-1}$ with the same minimal free resolution as that of $X_\L$. Let $T = k[x_0, \dots, x_{N-1}]$ be the homogeneous coordinate ring of $\pp^{N-1}$, and let $A$ be the homogeneous coordinate ring of $C$ in $\pp^{N-1}$. Then, the Betti number of $X_\L$ (also of $C$) are given by
$$\beta_{i,j} = \tor_i^T(A, k)_j, \ \text{for all} \ i,j.$$

Let $\H_X$ be the Hilbert function of $X_\L$ in $\pp^N$, i.e. 
$$\H_X(h) = \dim_k (S/I)_h.$$
Since $S/I$ is a Cohen-Macaulay ring, we have
$$\dim_k (S/I)_h = h^0(X_\L, \O_{X_\L}(h)) = h^0(X, \L^{\otimes h}).$$ 
As before, by Kodaira's vanishing theorem, $H^1(X, \L^{\otimes h}) = H^2(X, \L^{\otimes h}) = 0$ for $h > 0$. Thus, by the Riemann-Roch theorem, we get
\begin{align}
\dim_k (S/I)_h = h^0(X, \L^{\otimes h}) = \sfrac{1}{2}\L^{\otimes h}.(\L^{\otimes h} - K) + 1 + p_a \ \forall \ h > 0. \label{hilbertfn}
\end{align}
Let $d$ be the degree of $X_\L$ in $\pp^N$, then by (\ref{hilbertfn}),
\begin{align}
d \le \L^2 = (K+L)^2 = K^2 + 2K.L + L^2. \label{deg}
\end{align}

Let $D = C \cap H$ be a general hyperplane section of $C$. Then, $D$ is a set of $d$ points in $\pp^{N-2}$. Let $\J_C$ be the ideal sheaf of $C$ in $\pp^{N-1}$, and let $\J_D$ be the ideal sheaf of $D$ in $H \simeq \pp^{N-2}$. Let $\H_D$ be the Hilbert function of $D$. Since $C$ is arithmetically Cohen-Macaulay, it follows from the exact sequence
$$0 \To \J_C \To \J_C(1) \To \J_D(1) \To 0$$
that
$$0 \To H^1(\pp^{N-2}, \J_D(1)) \To H^2(\pp^{N-1}, \J_C) \To H^2(\pp^{N-1}, \J_C(1)) \To 0.$$
This implies that
\begin{align}
h^2(\pp^{N-1}, \J_C) - h^2(\pp^{N-1}, \J_C(1)) = h^1(\pp^{N-2}, \J_D(1)) = d - \H_D(1). \label{cohomology}
\end{align}

\begin{claim} $h^2(\pp^{N-1}, \J_C(1)) = 0$.
\end{claim}

\noindent {\it Proof of Claim.} Suppose $h^2(\pp^{N-1}, \J_C(1)) \not= 0$. By Kodaira's vanishing theorem, we have $H^1(X, \L) = H^2(X, \L) = 0$. Thus, by the Riemann-Roch theorem, we get
$$N = h^0(X, \L) - 1 = \sfrac{1}{2}\L.(\L-K) + p_a.$$
This implies that 
$$2N - 3 = \L(\L-K) + 2p_a - 3 = (K+L).L + 2p_a - 3 = L^2 + K.L + 2p_a-3.$$
Together with (\ref{deg}) we have
$$2N - 3 - d = - K.L - K^2 + 2p_a - 3 \ge p, $$
i.e.
$2N - 3 \ge d+p.$
It now follows from \cite{ba} that
$$\H_D(1) \ge \min \{ d, N-1 \} \ge d - (N-2) + p.$$
This and (\ref{cohomology}) give us
$$h^2(\pp^{N-1}, \J_C) - h^2(\pp^{N-1}, \J_C(1)) \le (N-2) - p < N-2.$$
Thus,
$$h^1(\pp^{N-1}, \O_C) - h^2(\pp^{N-1}, \O_C(1)) < N-2.$$
By Serre-Grothendieck correspondence, we get
\begin{align}
\dim_k \big[H^2_\m(A)\big]_0 - \dim_k \big[H^2_\m(A)\big]_1 < N-2, \label{contradiction}
\end{align}
where $\m$ is the maximal homogeneous ideal of $A$.

Let $x$ be a non-zero linear form in $A$. Since $A$ is a domain, we have the following exact sequence
$$0 \To A(-1) \stackrel{\times x}{\To} A \To A/xA \To 0.$$
This, by taking the corresponding long exact sequence of local cohomology and since $\dim A/xA = \dim A - 1 = 1$, gives an exact sequence
$$H^2_\m(A)(-1) \stackrel{\times x}{\To} H^2_\m(A) \To 0.$$
This and a special case of \cite[Lemma 3.1]{brod1} (cf. \cite[Lemma 3]{brod2}) imply that
$$\dim_k \big[H^2_\m(A)\big]_1 \le \dim_k \big[H^2_\m(A)\big]_0 - (N-2).$$
We obtain a contradiction to (\ref{contradiction}). The claim is proved. \qed

Now, we are back to the proof of our theorem. Having $h^2(\pp^{N-1}, \J_C(1)) = 0$, (\ref{cohomology}) gives us
$$h^2(\pp^{N-1}, \J_C) = h^1(\pp^{N-2}, \J_D(1)) = d - \H_D(1) \le (N-2) - p.$$
Let $\omega_A$ be the canonical module of $A$. Then,
$$\dim_k (\omega_A)_0 = \dim_k \big[H^2_\m(A)]_0 = h^1(\pp^{N-1}, \O_C) = h^2(\pp^{N-1}, \J_C) \le (N-2) - p.$$
By Lemma \ref{can-mod}, $\omega_A$ is a torsion free $A$-module. Thus, by a vanishing theorem of Eisenbud and Koh \cite[Theorem 1.1]{ek} (and its proof), $\omega_A$ has initial degree at least 0 and we have
$\big[\tor_s^T(\omega_A, k)\big]_s = 0 \ \text{for all} \ s \ge (N-2) - p.$ Therefore, 
$$\big[\tor_i^T(\omega_A, k)\big]_j = 0 \ \forall \ j < i \ \text{and} \ \big[\tor_s^T(\omega_A, k)\big]_s = 0 \ \forall \ s \ge (N-2) - p.$$
By duality, we now obtain
$$\reg \J \le 3 \ \text{and} \ \big[\tor_i^T(A, k)\big]_{i+2} = 0 \ \text{for} \ i \le p.$$
Hence, 
$$\beta_{i,j} = 0 \ \text{for} \ i \le p \ \text{and} \ j \ge i+2.$$
This implies that $\L$ has property $N_p$. The theorem is proved.
\end{proof}

\begin{cor} \label{surface-main}
Let $X$ be a surface with $H^1(X, \O_X) = 0$. Suppose $K$ is the canonical divisor and $L$ is an ample divisor on $X$ such that $\L = K \otimes L$ is generated by global sections. Then, for any integer $p \ge 0$, the adjoint line bundle $\L$ has property if $-K.L \ge p+3+K^2$.
\end{cor}

\begin{proof} We only need to show that $\L$ is normally generated; the conclusion then follows from Theorem \ref{surface} since $p_a \ge 0$. We need to prove that the natural maps
$$H^0(X, \L^{\otimes r}) \otimes H^0(X, \L) \To H^0(X, \L^{\otimes (r+1)})$$
are surjective for all $r \ge 1$. 

As before, by Kodaira's vanishing theorem and from the hypotheses, we have 
$$H^1(X, \L^{\otimes (r-1)}) = 0 \ \text{for all} \ r \ge 1.$$
Let $C \in |\L|$ be a smooth and irreducible curve. By \cite[Observation 1.1]{gp6}, it suffices to prove that the maps
$$H^0(C, \L^{\otimes r} \otimes \O_C) \otimes H^0(C, \L \otimes \O_C) \To H^0(C, \L^{\otimes (r+1)} \otimes \O_C)$$
are surjective for all $r \ge 1$. Since $-K.L \ge 3+K^2$, we have $-K.\L \ge 3$, and so by adjunction formula, $\deg (\L \otimes \O_C) \ge 2g(C) + 1$. Hence, the assertion follows from Castelnuovo's theorem.
\end{proof}

As a consequence of Corollary \ref{surface-main}, we can now prove Mukai's conjecture for rational anticanonical surfaces. Our argument is a combination of that of \cite[Proposition 1.10]{gp6} and the following lemma.

\begin{lemma} \label{normal}
Let $X$ be a regular anticanonical surface. Let $K$ and $L$ be the canonical divisor and an ample line bundle on $X$. Then
\begin{enumerate}
\item $L$ is generated by global sections if $-K.L \ge 2$,
\item $L$ is normally generated if $-K.L \ge 3$.
\end{enumerate}
\end{lemma}

\begin{proof} By the assumption, $|-K|$ has an effective divisor. This implies that $-K.L > 0$. Thus, the Kodaira dimension of $X$ must be negative. It then follows from \cite[Theorem V.6.1]{hart} that $X$ is a rational or a ruled surface. In the later case, since $H^1(X, \O_X) = 0$, $X$ must be a rational ruled surface. (1) now follows from \cite[Theorem III.1]{harb}.

To prove (2), we first observe that for any $r \ge 1$, $L^{\otimes r}$ is globally generated by (1). Thus, there exists a smooth irreducible curve $C_r \in |L^{\otimes r}|$. We have an exact sequence
$$0 \rightarrow \O_X \rightarrow L^{\otimes r} \rightarrow L^{\otimes r} \otimes \O_{C_r} \rightarrow 0.$$
Taking the associated long exact sequence of cohomology groups, we get
\begin{align}
0 \rightarrow H^1(X, L^{\otimes r}) \rightarrow H^1(C_r, L^{\otimes r} \otimes \O_{C_r}). \label{eq.normal1}
\end{align}
Since $K.L^{\otimes r} = r(K.L) < 0$, by adjunction formula, $\deg (L^{\otimes r} \otimes \O_{C_r}) > 2g(C_r) - 2$. Thus, $H^1(C_r, L^{\otimes r} \otimes \O_{C_r}) = 0$. This and (\ref{eq.normal1}) imply that $H^1(X, L^{\otimes r}) = 0$. This holds for any $r \ge 1$. Together with the hypothesis that $X$ is regular, we have
$$H^1(X, L^{\otimes (r-1)}) = 0, \ \forall \ r \ge 1.$$

Now, by (1), $L$ is generated by global sections. We, therefore, can make use of \cite[Observation 1.1]{gp6} in the same line of arguments as that of Corollary \ref{surface-main} to conclude that $L$ is normally generated.
\end{proof}

\begin{thm} \label{surface-np}
Let $X$ be a rational anticanonical surface. Let $K$ be the canonical divisor and let $L_1, \dots, L_m$ be ample line bundles on $X$. Then, for any integer $p \ge 0$, the adjoint line bundle $K \otimes L_1 \otimes \dots \otimes L_m$ has property $N_p$ if $n \ge 4+p$.
\end{thm}

\begin{proof} Let $L = L_1 \otimes \dots \otimes L_m$. Since $X$ is anticanonical, $-K$ is in the same class with an effective divisor. Thus, $-K.L_i \ge 1$ for all $i = 1, \dots, m$. Without loss of generality, suppose that 
$$s = -K.L_1 = \min \{ -K.L_i ~|~ i=1, \dots, m \}.$$
Then, $s \ge 1$ and
$$-K.L = \sum_{i=1}^m -K.L_i \ge ms \ge (4+p)s.$$
It follows from Reider's theorem \cite{reider} that $K \otimes L$ is very ample. In particular, $K \otimes L$ is generated by global sections. By Corollary \ref{surface-main}, it suffices to show that for any integer $p \ge 0$,
\begin{align}
(4+p)s \ge p+3+K^2. \label{inequality}
\end{align}

If $-K.(K \otimes L_1) \ge 0$ then $s \ge K^2$, and so $(4+p)s = (3+p)s + s \ge p+3+K^2$, i.e. (\ref{inequality}) holds. Assume that $-K.(K \otimes L_1) < 0$. By Kodaira's vanishing theorem, we have $H^1(X, K \otimes L_1) = H^2(X, K \otimes L_1) = 0$. Thus, the Riemann-Roch theorem gives us
\begin{align}
h^0(X, K \otimes L_1) = \sfrac{1}{2}(K \otimes L_1).L_1 + 1 + p_a. \label{h0}
\end{align}
Suppose $(K \otimes L_1).L_1 \ge 0$. By (\ref{h0}), $h^0(X, K \otimes L_1) > 0$, i.e. $K \otimes L_1$ is in the same class with an effective divisor. Thus, $-K.(K+L_1) \ge 0$, and we are done. From now on, we shall assume that 
\begin{align}
-K.(K \otimes L_1) < 0 \ \text{and} \ (K \otimes L_1).L_1 < 0. \label{negative}
\end{align}

If $K^2 \le 1$, then (\ref{inequality}) holds since $s \ge 1$. Suppose now that $K^2 \ge 2$. Since $L_1$ is ample, $L_1^2 > 0$. We will consider different cases based on the value of $L_1^2$.

Assume first that $L_1^2 \ge 5$. By Hodge Index inequality, we must have $-K.L_1 > 3$. By Lemma \ref{normal}, $L_1$ is normally generated. This, in particular, implies that $L_1$ is very ample. By Reider's theorem \cite{reider}, $K \otimes L_1$ is generated by global sections unless there exists a reduced irreducible curve $E$ such that $L_1.E = 1$ and $E^2 = 0$. If $K \otimes L_1$ is generated by global sections, then $h^0(X, K \otimes L_1) > 0$, i.e. $K \otimes L_1$ is in the same class with an effective divisor, and so we have $-K.(K+L_1) \ge 0$. This is a contradiction to the assumption (\ref{negative}). If there exists a reduced and irreducible curve $E$ such that $L_1.E = 1$ and $E^2 = 0$, then since $L_1$ is very ample, $E = \pp^1$. It follows from \cite[Lemma 4.1.10]{bs} that $X$ is a $\pp^1$-bundle over a curve of genus 0. This implies that $K^2 = 8$. We now have $(4+p)s > (4+p)3 \ge p+12 > p+3+K^2$. (\ref{inequality}) holds.

Assume that $L_1^2 = 4$. By Hodge Index inequality again, we have $-K.L_1 \ge 3$. As before, by Lemma \ref{normal}, this implies that $L_1$ is normally generated. Thus, $L_1$ embeds $X$ as a projectively normal surface. However, the only normal smooth surfaces of degree 4 are K3 surfaces in $\pp^3$, the Del Pezzo surface in $\pp^4$, and the Veronese surface and the normal scroll in $\pp^5$ (cf. \cite[pp. 145, 173]{io}). If $X$ is embedded as a K3 surface then $K^2 = 0$, a contradiction to the assumption that $K^2 \ge 2$. If $X$ is embedded as a Del Pezzo surface by $L_1$ then $K = -L_1$, and so $-K(K \otimes L_1) = 0$, a contradiction to (\ref{negative}). If $X$ is embedded as a Veronese surface then $K^2 = 9$, and we establish (\ref{inequality}) since $(4+p)s \ge (4+p)3 \ge p+12 = p+3+K^2$. If $X$ is embedded as a normal scroll in $\pp^5$ then $K^2 = 8$, and again we establish (\ref{inequality}). 

Assume that $L_1^2 = 3$. By Hodge Index inequality, we again have $-K.L_1 \ge 3$, and so $L_1$ is normally generated. This implies that $X$ is either a Del Pezzo cubic surface in $\pp^3$ or a normal scroll in $\pp^4$ (cf. \cite[pp. 145, 173]{io}). As before, we either arrive with a contradiction to (\ref{negative}) or establish (\ref{inequality}).

Assume that $L_1^2 = 2$. By Hodge Index inequality, we have $-K.L_1 \ge 2$. If $-K.L_1 = 2$ then $(K \otimes L_1).L_1 = 0$, a contradiction to (\ref{negative}). Suppose $-K.L_1 \ge 3$. By Lemma \ref{normal} again, $L_1$ is normally generated. In particular, $L_1$ is very ample. This implies that $(X, L_1) = (\pp^1 \times \pp^1, \O_{\pp^1 \times \pp^1}(1))$ (cf. \cite[pp. 145, 173]{io}). Thus, $K^2 = 8$ and we establish (\ref{inequality}). 

Assume finally that $L_1^2 = 1$. As before, by Kodaira's vanishing theorem and the Riemann-Roch theorem, we have 
$h^0(X, K \otimes L_1) = \sfrac{1}{2}L_1.(K \otimes L_1) + 1 + p_a.$
Thus, $L_1.(K \otimes L_1)$ is an even number, i.e. $K.L_1$ is odd. By Hodge Index inequality, $-K.L_1 \ge 2$. Thus, $-K.L_1 \ge 3$. This, once again, implies that $L_1$ is very ample. Since $L_1^2 = 1$, we must have $(X, L_1) = (\pp^2, \O_{\pp^2}(1))$ (cf. \cite[pp. 145, 173]{io}). In this case $K^2 = 9$, and we establish (\ref{inequality}).

The theorem is proved.
\end{proof}

As a consequence, Theorem \ref{surface-np} establishes Mukai's conjecture for rational anticanonical surfaces.

\begin{cor} \label{surface-cor}
Let $X$ be a rational anticanonical surface. Let $K$ be the canonical divisor and let $L$ be an ample line bundle on $X$. Then, for any integer $p \ge 0$, the adjoint line bundle $K \otimes L^{\otimes (4+p)}$ has property $N_p$.
\end{cor}

\begin{proof} The assertion follows from Theorem \ref{surface-np} by taking $L_1 = \dots = L_m = L$.
\end{proof}


\section{Regularity of ideal sheaves of embeddings} \label{section-reg}

\noindent In this section, we study width of the minimal free resolution of embeddings of a regular projective variety given by adjoint line bundles. This is done via examining the regularity of ideal sheaves of these embeddings. The reader is again referred to Definition \ref{aCM-def} for the notion of regular projective varieties. 

Our main result in this section is stated in the next theorem. Parts (1) and (2) follow straight forward from a standard use of Kodaira's vanishing theorem. Part (3) is what we are interested in.

\begin{thm} \label{aCM}
Let $X$ be a smooth projective variety of dimension $n$. Suppose $K$ is the canonical divisor and $L$ is an ample line bundle on $X$ such that $\L = K \otimes L$ is normally generated. Let $\J_\L$ be the ideal sheaf of the embedding of $X$ given by $\L$. Then,
\begin{enumerate}
\item $\reg \J_\L \le n+2.$ Moreover, if $H^0(X, K) \not= 0$, then $\reg \J_\L = n+2$.
\item If $H^0(X, K) = 0$, then $\reg \J_\L \le n+1$.
\item If $X$ is regular, $H^0(X, K) = 0$ and $H^0(X, K \otimes \L) \not= 0$, then
$\reg \J_\L = n+1.$
\end{enumerate}
\end{thm}

\begin{proof} (1) Since $\L$ is normally generated, $\L$ is necessarily very ample (cf. \cite{mum}). Let $N = \dim_k H^0(X, \L) - 1$, and let $X_\L$ be the image of $X$ embedded in $\pp^N$ by $\L$. Let $S$ and $I$ be the coordinate ring and the defining ideal of $X_\L$ in $\pp^N$. We have the following exact sequence
$$0 \To I \To S \To \oplus_{h \ge 0} H^0(X, \L^{\otimes h}) \To \oplus_{h \ge 0} H^1(\pp^N, \J_\L(h)) \To 0.$$
Since $\L$ is normally generated, $\L$ embeds $X$ as a projectively normal variety. Thus, we have
\begin{align}
S/I \simeq \oplus_{h \ge 0} H^0(X, \L^{\otimes h}) \ \text{and} \ H^1(\pp^N, \J_\L(h)) = 0 \ \forall \ h \ge 0. \label{1}
\end{align}

By taking the long exact sequence of sheaf cohomology associated to the exact sequence 
$0 \To \J_\L \To \pp^N \To \O_{X_\L} \To 0,$
we obtain
\begin{align}
H^i(\pp^N, \J_\L(h)) = H^{i-1}(\pp^N, \O_{X_\L}(h)) \ \text{for all} \ i \ge 2 \ \text{and} \ h \ge -N. \label{structure}
\end{align}
For $h > 0$, we have 
$\L^{\otimes h} = K \otimes \L^{\otimes (h-1)} \otimes L$.
Since $L$ and $\L$ are ample, by Kodaira's vanishing theorem, we get
\begin{align}
H^i(X, \L^{\otimes h}) = 0, \ \forall \ i \ge 1, \ \text{and} \ h > 0. \label{kodaira}
\end{align}
Thus, $H^i(\pp^N, \O_{X_\L}(n+1-i)) = H^i(X, \L^{\otimes (n+1-i)}) = 0$ for $1 \le i \le n$. Clearly, $H^i(\pp^N, \O_{X_\L}(n+1-i)) = 0$ for $i > n = \dim X_\L$. Hence, (\ref{structure}) gives us
$$H^i(\pp^N, \J_\L(n+2-i)) = 0 \ \text{for all} \ 2 \le i \le N.$$
This and (\ref{1}) prove that 
$$\reg \J_\L \le n+2.$$

Now, if $H^0(X, K) \not= 0$, then by Serre's duality, we have $H^n(X, \O_X) = H^0(X, K) \not= 0$. Thus, by (\ref{structure}), 
$$H^{n+1}(\pp^N, \J_\L) = H^n(\pp^N, \O_{X_\L}) = H^n(X, \O_X) \not= 0.$$
This implies that $\reg \J_\L \ge n+2$. Hence, $\reg \J_\L = n+2$.

\noindent (2) By Serre's duality, we have $H^n(\pp^N, \O_{X_\L}) = H^n(X, \O_X) = H^0(X, K) = 0$. This and (\ref{kodaira}) imply that
$H^i(\pp^N, \O_{X_\L}(n-i)) = H^i(X, \L^{\otimes (n-i)}) = 0 \ \text{for all} \ 1 \le i \le n.$ It's again clear that $H^i(\pp^N, \O_{X_\L}(n-i)) = 0$ for $i > n = \dim X_\L$.
It now follows from (\ref{structure}) that
\begin{align}
H^i(\pp^N, \J_\L(n+1-i)) = 0 \ \text{for all} \ 2 \le i \le N. \label{2}
\end{align}
The conclusion follows from (\ref{1}) and (\ref{2}).

\noindent (3) Suppose that $X$ is regular, $H^0(X, K) = 0$ and $H^0(X, K \otimes \L) \not= 0$. By (2), $\reg \J_\L \le n+1$. We need to show that
\begin{align}
\reg \J_\L \ge n+1. \label{n+1}
\end{align}
By (\ref{1}), $S/I \simeq \oplus_{h \ge 0} H^0(X, \L^{\otimes h})$. Thus, by Green's Syzygy theorem \cite[1.b.4]{gr}, $S/I$ admits a minimal free resolution
$$0 \To F_s \To \dots \To F_1 \To F_0 = S \To S/I \To 0, $$
for some $s \ge N-n$ (the codimension of $X_\L$ in $\pp^N$) and 
$$F_i = \bigoplus_{q \ge 1} \K_{i,q}(\L) \otimes S(-i-q), \ \text{for} \ i =1, \dots, s,$$
where $\K_{i,q}(\L)$ is the Koszul cohomology group associated to $\L$ as in Definition \ref{koszul}.
It suffices to show that there exists an integer $p \ge 0$ such that $$\K_{p, n}(\L) \not= 0.$$

As in (\ref{kodaira}), we have 
$H^i(X, \L^{\otimes h}) = 0 \ \text{for} \ i \ge 1 \ \text{and} \ h > 0.$
This, together with the hypothesis that $X$ is regular, gives us
$$H^i(X, \L^{\otimes (n-i)}) = H^i(X, \L^{\otimes (n-i-1)}) = 0 \ \text{for} \ i = 1, \dots, n-1.$$
Therefore, Green's Duality theorem \cite[2.c.6]{gr} applies for $p \ge 0$ and $q = n$, and gives
$$\K_{N-n, n}(\L)^* \simeq \K_{0,1}(K, \L).$$
We shall prove that $\K_{0,1}(K, \L) \not= 0$ under the given hypotheses. Indeed, let $W = H^0(X, \L)$. The Koszul complex of $\oplus_{h \ge 0} H^0(X, K \otimes \L^{\otimes h})$ at degree $(0,1)$ is
$$ \bigwedge W \otimes H^0(X, K) \stackrel{d_{1,0}}{\To} W \otimes H^0(X, K \otimes \L) \stackrel{d_{0,1}}{\To} 0.$$
Hence, since $H^0(X, K) = 0$ and $H^0(X, K \otimes \L) \not= 0$, 
$$\K_{0,1}(K, \L) = W \otimes H^0(X, K \otimes \L) \not= 0.$$
The theorem is proved. 
\end{proof}

Theorem \ref{aCM}, combined with Fujita's freeness and very ampleness conjectures (cf. \cite{f, ko2, smith}), gives rise to several interesting corollaries. It is well known that Fujita's freeness and very ampleness conjectures hold for ample line bundles which are generated by global sections. For completeness, we shall still include the proof in the next proposition.

\begin{pro} \label{fujita}
Let $X$ be a smooth projective variety of dimension $n$ with canonical divisor $K$. Let $L$ be an ample line bundle which is generated by global sections. Then, $K \otimes L^{\otimes m}$ is generated by global sections for $m \ge n+1$, and is very ample for $m \ge n+2$.
\end{pro}

\begin{proof} Let $\L_m = K \otimes L^{\otimes m}$ for $m \ge 0$. For $h \ge n+1$, we have $\L_h \otimes L^{\otimes (-i)} = K \otimes L^{\otimes (h-i)}$. It follows from Kodaira's vanishing theorem that
$$H^i(X, \L_h \otimes L^{\otimes (-i)}) = 0$$
for $i = 1, \dots, n$. Thus, $\L_h$ is $0$-regular in the sense of Castelnuovo-Mumford with respect to $L$. Since $L$ is generated by global sections, this implies (cf. \cite[Theorem 1.8.5 and Example 1.8.22]{l}) that $\L_h$ is generated by global sections and $\L_h \otimes L$ is very ample. In other words, $\L_m$ is generated by global sections for $m \ge n+1$, and is very ample for $m \ge n+2$.
\end{proof}

\begin{cor} \label{veryample}
Let $X$ be a smooth regular projective variety of dimension $n$. Let $K$ be the canonical divisor and let $L$ be a very ample line bundle on $X$. Let $\L_m = K \otimes L^{\otimes m}$, and let $\J_m$ be the ideal sheaf of the embedding of $X$ given by $\L_m$ for $m \ge n+2$. The following statements hold.
\begin{enumerate}
\item If $H^0(X, K) \not= 0$, then $\reg \J_m = n+2$ for $m \ge n+2$.
\item If $H^0(X, K) = 0$, then
$\reg \J_m = n+1 \ \text{for all} \ m \ge 2n+2.$
\end{enumerate}
\end{cor}

\begin{proof} By \cite[Theorem 1]{el}, $\L_m$ is normally generated for $m \ge n+2$. The first conclusion follows from Theorem \ref{aCM}. Observe further that by Proposition \ref{fujita} for $m \ge 2n+2$, $\L_{n+1}$ and $\L_{m-n-1}$ are generated by global sections, and so $K \otimes \L_m = \L_{n+1} \otimes \L_{m-n-1}$ is generated by global sections. Thus, $H^0(K \otimes \L_m) \not= 0$ for $m \ge 2n+2$. The second assertion is also a consequence of Theorem \ref{aCM}.
\end{proof}

\begin{remark} \label{behaviour}
It is worth noting that having the same hypotheses as in Corollary \ref{veryample}, by \cite[Theorem 1]{el} we know that for $m \ge 2n+2$, $\L_m$ has property $N_{m-n-2}$. This and Corollary \ref{veryample} say that as $m$ gets large, while the beginning part of the minimal free resolution of $\J_m$ becomes more linear, the width of the resolution is always fixed, either $n+2$ or $n+1$. 
\end{remark}

\begin{cor} \label{pn}
Let $X$ be a smooth rational projective variety of dimension $n$, and let $K$ be the canonical divisor of $X$. Suppose $L$ is a very ample line bundle on $X$. Let $\L_m = K \otimes L^{\otimes m}$, and let $\J_m$ be the ideal sheaf of the embedding of $X$ given by $\L_m$ for $m \ge n+2$. Then,
$$\reg \J_m = n+1 \ \text{for} \ m \ge 2n+2.$$
\end{cor}

\begin{proof} Since $X$ is rational, $X$ is obtained from $\pp^n$ by blowing up along an ideal sheaf. Let $\pi: X \To \pp^n$ be the blowing up along an ideal sheaf $\I$ associated to an ideal $I \subset k[x_0, \dots, x_n]$. Let $E_0$ and $E$ be the pull-back of a general hyperplane in $\pp^n$ and the exceptional divisor of $\pi$. Let $r$ be the height of $I$. We have
$$K = -(n+1)E_0 + (r-1)E.$$
Thus, $H^0(X, K) = 0$. Moreover, it follows from Remark \ref{rational} that
$X$ is regular.
Hence, the conclusion follows from Corollary \ref{veryample}.
\end{proof}

\begin{remark} For surfaces, making use of Reider's theorem \cite{reider} and results in a series of papers of Gallego and Purnaprajna \cite{gp1, gp2, gp3, gp4, gp5, gp6}, we can relax the condition that $L$ is very ample in Corollary \ref{veryample}. Yet, the bounds obtained would be slightly weaker since Mukai's conjecture for projective normality (property $N_0$) is not completely known to hold.
\end{remark}



\begin{thebibliography}{10}

\bibitem{as} Angehrn, U. and Siu, Y.T. {\it Effective freeness and point separation for adjoint bundles}. Invent. Math. {\bf 122} (1995), 291-308.
\bibitem{ba} Ballico, E. {\it On singular curves in positive characteristic}. Math. Nachr. {\bf 141} (1989), 267-273.
\bibitem{bs} Beltramitti, M.C. and Sommese, A.J. The adjunction theory of complex projective varieties. Walter de Gruyter, 1995.
\bibitem{brod1} Brodmann, M. {\it Bounds on the cohomological Hilbert functions of a projective variety}. J. Algebra. {\bf 109} (1987), 352-380.
\bibitem{brod2} Brodmann, M. {\it A bound for the first cohomology of a projective surface}. Arch. Math. {\bf 50} (1988), 68-72.
\bibitem{chtv} Conca, A., Herzog, J., Trung, N.V., and Valla, G. {\it Diagonal subalgebras of bigraded algebras and embeddings of blow-ups of projective spaces}. Amer. J. Maths. {\bf 119} (1997), 859-901.
\bibitem{cha} Cutkosky, S.D. and H\`a, H. T\`ai. {\it Arithmetic Macaulayfication of projective schemes}. J. Pure Appl. Algebra. {\bf 201} (2005), 49-61.
\bibitem{che} Cutkosky, S.D. and Herzog, J. {\it Cohen-Macaulay coordinate rings of blowup schemes}. Comment. Math. Helv. {\bf 72} (1997), 605-617.
\bibitem{d} Demailly, J.P. {\it A numerical criterion for very ample line bundles}. J. Diff. Geom. {\bf 37} (1993), 323-374.
\bibitem{de} Demazure, M. {\it Surfaces de Del Pezzo - IV. Syst\`emes anticanoniques}. S\'eminaire sur les Singularit\'es des surfaces. Lecture Notes in Mathematics. {\bf 777} (1980).
\bibitem{el} Ein, L. and Lazarsfeld, R. {\it Syzygies and Koszul cohomology of smooth projective varieties of arbitrary dimension}. Invent. Math. {\bf 111} (1993), 51-67.
\bibitem{ek} Eisenbud, D. and Koh, J. {\it Some linear syzygy conjectures}. Adv. Math. {\bf 90} (1991), 47-76.
\bibitem{f} Fujita, T. {\it Defining equations for certain types of polarized variety}. In Complex Analysis and Algebraic Geometry, W.L. Baily, Jr. and T. Shioda, eds. 1977. Cambridge Univ. Press. 165-173.
\bibitem{gp1} Gallego, F.J., Purnaprajna, B.P. {\it Normal presentation on elliptic ruled surfaces}. J. Algebra. {\bf 186} (1996), no. 2, 597-625.
\bibitem{gp2} Gallego, F.J., Purnaprajna, B.P. {\it Higher syzygies of elliptic ruled surfaces}. J. Algebra. {\bf 186} (1996), no. 2, 626-659.
\bibitem{gp3} Gallego, F.J. and Purnaprajna, B.P. {\it Projective normality and syzygies of algebraic surfaces}. J. Reine Angew. Math. {\bf 506} (1999), 145-180.
\bibitem{gp4} Gallego, F.J., Purnaprajna, B.P. {\it Syzygies of projective surfaces: an overview}. J. Ramanujan Math. Soc. {\bf 14} (1999), no. 1, 65-93.
\bibitem{gp5} Gallego, F.J., Purnaprajna, B.P. {\it Vanishing theorems and syzygies for $K3$ surfaces and Fano varieties}. J. Pure Appl. Algebra. {\bf 146} (2000), no. 3, 251-265.
\bibitem{gp6} Gallego, F.J., Purnaprajna, B.P. {\it Some results on rational surfaces and Fano varieties}. J. Reine Angew. Math. {\bf 538} (2001), 25-55.
\bibitem{ggp} Geramita, A.V., Gimigliano, A. and Pitteloud, Y. {\it Graded Betti
numbers of some embedded rational $n$-folds}. Math. Ann. {\bf 301} (1995), 363-380.
\bibitem{gr} Green, M. (1984). {\it Koszul cohomology and the geometry of projective varieties I, II}. J. Diff. Geom. {\bf 19, 20}, 125-171, 279-289.
\bibitem{gl} Green, M. and Lazarsfeld, R. {\it Some results on the sygyzies of finite sets and algebraic curves}. Comp. Math. {\bf 67} (1988), 301-314.
\bibitem{ha1} H\`a, Huy T\`ai. {\it Projective embeddings of projective schemes blown up at subschemes}. Math. Z. {\bf 246} (2004), no. 1-2, 111-124.
\bibitem{ht} H\`a, Huy T\`ai and Trung, N.V. {\it Asymptotic behaviour of arithmetically Cohen-Macaulay blow-ups}. Trans. Amer. Math. Soc. {\bf 357} (2005), no. 9, 3655-3672.
\bibitem{harb} Harbourne, B. {\it Anticanonical rational surfaces}. Trans. Amer. Math. Soc. {\bf 349} (1997), no. 3, 1191-1208.
\bibitem{hart} Hartshorne, R. Algebraic Geometry. Graduate Text {\bf 52}. Springer-Verlag, 1977.
\bibitem{ss} Hering, M., Schenck, H., and Smith, G.G. {\it Syzygies, multigraded regularity and toric varieties}. Compositio Math. {\bf 142} (2006), 1499-1506.
\bibitem{io} Ionescu, P. {\it Embedded projecitve varieties of small invariants}. Lecture Notes in Mathematics. {\bf 1056} (1982), 142-186.
\bibitem{ko2} Koll\'ar, J. {\it Singularities of pairs}. Proc. of Symposia in Pure Math. {\bf 62} (1997), no. 1, 221-287.
\bibitem{la1} Lazarsfeld, R. {\it A sampling of vector bundles techniques in the study of linear series.} In Cornalba, M., Gomez-Mont, X., Verjovsky, A. (eds): Lectures on Riemann surfaces, World Scientific (1989), 500-559.
\bibitem{l} Lazarsfeld, R. Positivity in algebraic geometry. Book to appear.
\bibitem{mum} Mumford, D. {\it Varieties defined by quadratic equations}. In Questions on algebraic varieties. Proc. CIME 1970, ed. Cremonese, Roma, 31-100.
\bibitem{pa} Pareschi, G. {\it Syzygies of abelian varieties.} J. Amer. Math. Soc. {\bf 13} (2000), no. 3, 651-664. 
\bibitem{pp1} Pareschi, G. and Popa, M. {\it Regularity on abelian varieties. I.} J. Amer. Math. Soc. {\bf 16} (2003), no. 2, 285-302.
\bibitem{pp2} Pareschi, G. and Popa, M. {\it Regularity on abelian varieties. II. Basic results on linear series and defining equations.} J. Algebraic Geom. {\bf 13} (2004), no. 1, 167-193.  
\bibitem{p1} Park, E. {\it On higher syzygies of ruled surfaces}. Preprint. aXiv:math.AG/0401100.
\bibitem{p2} Park, E. {\it On higher syzygies of ruled surfaces II}. Preprint. aXiv:math.AC/0411022.
\bibitem{reider} Reider, I. {\it Vector bundles of rank 2 and linear systems on algebraic surfaces}. Ann. Math. {\bf 127} (1988), 309-316.
\bibitem{rubei} Rubei, E. {\it On syzygies of abelian varieties.} Trans. Amer. Math. Soc. {\bf 352} (2000), no. 6, 2569-2579.
\bibitem{smith} Smith, K.E. {\it A tight closure proof of Fujita's freeness conjecture for very ample line bundles}. Math. Ann. {\bf 317} (2000), 285-293.
\bibitem{std} St. Donat, B. {\it Sur les \'equations d\'efinissant une courbe alg\'ebrique}. C. R. Acad. Sci. Paris. {\bf 274} (1972), 324-327.

\end{thebibliography}
\end{document}